\documentstyle[12pt]{article}
\if@twoside  m
\else
\fi
\newcommand{\ie}{{\em i.e.}}
\newcommand{\eg}{{\em e.g.}}
\newcommand{\cf}{{\em cf. }}
\newcommand{\QED}{\mbox{\rule[-1.5pt]{6pt}{10pt}}}
\newcommand{\rhs}{{\em rhs }}
\newcommand{\R}{I\!\!R}
\newcommand{\HH}{{\cal H}}
\newcommand{\NN}{{\cal N}}
\newcommand{\eps}{\varepsilon}
\newtheorem{claim}{Claim}[section]
\newtheorem{theorem}[claim]{Theorem}
\newtheorem{proposition}[claim]{Proposition}
\newtheorem{remarks}[claim]{Remarks}

\begin{document}

\title{\Large\bf On the number of particles which a curved quantum
waveguide can bind}

\author{Pavel Exner$^{a,b}$ and Simeon A. Vugalter$^a$}
\date{}

\maketitle

\begin{center}
a) Nuclear Physics Institute, Academy of Sciences,
25068 \v{R}e\v{z} near Prague, \\
b) Doppler Institute, Czech Technical University, B\v rehov\'a 7,
11519 Prague, \\
{\em exner@ujf.cas.cz, vugalter@ujf.cas.cz} 
\end{center}

\begin{abstract}
We discuss the discrete spectrum of $\,N\,$ particles in a curved
planar waveguide. If they are neutral fermions, the maximum number of
particles which the waveguide can bind is given by a one--particle
Birman--Schwinger bound in combination with the Pauli principle. On
the other hand, if they are charged, \eg, electrons in a bent quantum
wire, the Coulomb repulsion plays a crucial role. We prove a
sufficient condition under which the discrete spectrum of such a
system is empty.
\end{abstract}


\section{Introduction}

A rapid progress of mesoscopic physics brought, in particular,
interesting new problems concerning relations between geometry and
spectral properties of quantum Hamiltonians. They involve models of
quantum wires, dots, and similar systems. While in reality these are
rather complicated systems composed of different semiconductor
materials, experience tells us that their basic features can be
explained using simple models in which electrons (regarded as free
particles with an effective mass) are supposed to be confined to an
appropriate spatial region, either by a potential or by a hard wall.
A brief description of this approximation with a guide to further
reading is given in Ref.~\cite{DE}. 
In addition, such models apply not only to electrons in semiconductor
microstructures; a different example is represented by atoms trapped
in hollow optical fibers \cite{SMZ}.

It is natural that most theoretical results up to date refer to the
case of a single particle in the confinement. On the other hand, from
the practical point of view it is rather an exception than a rule
that an experimentalist is able to isolate a single electron or atom,
and therefore many--body problems in this setting are of interest.
For instance, two--dimensional quantum dots which can be regarded as
artificial atoms have been studied recently, usually in presence of a
magnetic field, either for a pair of electrons or in the semiclasical
situation when a Thomas--Fermi--type approach is applicable --- \cf
[3-6]
and references therein.

In these studies, however, geometry of the dot played a little role,
because the confinement was realized by a harmonic potential or a
circular hard wall. This is not the case for open systems modelling
quantum wires where a deformation of a straight channel is needed to
produce nontrivial spectral properties. In particular, a quantum
waveguide exhibit bound states if it is bent \cite{DE,ES,GJ},
protruded 
[9-11]
or allowing a leak to another duct 
[12-14], 
and the
discrete spectrum depends substantially on the shape of the channel.
With few exceptions such as Ref.~\cite{NTV}, however, the known
results refer to the one--particle case.

It is the aim of the present paper to initiate a rigorous
investigation of many--particle effects in quantum waveguides.
We are going to discuss here a system of $\,N\,$ particles in
a bent planar Dirichlet tube, \ie, a hard--wall channel, and ask
whether $\,N$--particle bound states exist for a given geometry.
After collecting the necessary preliminaries in the next section,
we shall derive first in Section~3 a simple bound for the neutral
case which follows from the Birman--Schwinger estimate of the
one--particle Hamiltonian in combination with the Pauli principle.

The main result of the paper is formulated and proved in Section~4.
It concerns the physically interesting case of charged particles; the
example we have in mind is, of course, electrons in a bent
semiconductor quantum wire. The electrostatic repulsion makes
spectral analysis of the corresponding Hamiltonian considerably more
complicated. Using variational technique borrowed from atomic
physics, we derive here a sufficient condition under which the
discrete spectrum is empty. The condition is satisfied for $\,N\,$
large enough and represents an implicit equation for the maximum
number of charged particles which a waveguide of a given curvature
and width can bind. Some other aspects of the result and open
questions are discussed briefly in the concluding section.


\section{Preliminaries}

The waveguide in question will be modelled by a curved planar strip
$\,\Sigma\,$ in $\,\R^2\,$, of a constant width $\,d=2a\,$. It can be
obtained by transporting the perpendicular interval $\,[-a,a]\,$
along the curve $\,\Gamma\,$ which is the axis of $\,\Sigma\,$. Up
to Euclidean transformations, the strip is uniquely characterized by
its halfwidth $\,a\,$ and the (signed) curvature
$\,s\mapsto\gamma(s)\,$ of $\,\Gamma\,$, where $\,s\,$ denotes the
arc length. We adopt the regularity assumptions of
Refs.~\cite{DE,ES}: 
   \begin{description}
   \item{\em (i)} $\;\Omega\,$ is not self--intersecting,
   \item{\em (ii)} $\;a\|\gamma\|_{\infty}<1\,$,
   \item{\em (iii)} $\;\gamma\,$ is piecewise $\,C^2\,$ with
$\,\gamma',\, \gamma''\,$ bounded,
   \end{description}
and restrict our attention to the case when the tube is curved in a
bounded region only:
   \begin{description}
   \item{\em (iv)} there is $\,b>0\,$ such that $\,\gamma(s)=0\,$ for
$\,|s|>b\,$; without loss of generality we may assume that $\,2b>a\,$.
   \end{description}
As usual we put $\,\hbar=2m=1\,$; then the one--particle Hamiltonian
of such a waveguide is the Dirichlet Laplacian
$\,-\Delta_D^{\Sigma}\,$ defined in the conventional way --- \cf
\cite{RS}, Sec.~XIII.15. Using the natural locally orthogonal
curvilinear coordinates $\,s,u\,$ in $\,\Sigma\,$ one can map 
$\,-\Delta_D^{\Sigma}\,$ unitarily onto the operator 
\begin{equation} \label{Hamiltonian 1}
H_1\,=\,-\partial_s\, (1+u\gamma)^{-2}\, \partial_s\,-\,
\partial_u^2 \,+\,V(s,u) 
\end{equation}
on $\,L^2(\R\times(-a,a)\,)\,$ with the effective curvature--induced
potential 
\begin{equation} \label{effective potential}
V(s,u)\,:=\,-\,{\gamma(s)^2\over 4(1+u\gamma(s))^2}\,+\,
{u\gamma''(s) \over 2(1+u\gamma(s))^3}\,-\, 
{5 \over 4}\,{u\gamma'(s)^2 \over(1+u\gamma(s))^4} 
\end{equation}
which is e.s.a. on the core $\,D(H)\,=\, \{\,\psi\,:\,\psi\in
C^{\infty},\; \psi(s,\pm a)=0,\; H\psi\in L^2\,\}\,$ --- \cf
Refs.~\cite{DE,ES} for more details.

If the waveguide contains $\,N\,$ particles, the state Hilbert space
is $\,L^2(\Sigma))^N$; the Pauli principle will be taken into account
later. We assume that each particle has the charge $\,e\,$; using the
same ``straightening" transformation we are then able to rewrite the
Hamiltonian as
\begin{eqnarray} \label{Hamiltonian N}
H_N \,\equiv\, H_N(\gamma,a,e) &\!=\!&  \sum_{j=1}^N \left\lbrace\,
-\partial_{s_j}\, (1+u_j\gamma(s_j))^{-2}\, \partial_{s_j}\,-\, 
\partial_{u_j}^2 \,+\,V(s_j,u_j) \,\right\rbrace 
\nonumber \\ \nonumber \\
&\! +\!& e^2 \sum_{1\le j<l\le N} |\vec r_j- \vec r_l|^{-1}\,,
\end{eqnarray}
with the domain $\,\left(\HH^2(\R)\otimes \HH^2_0(-a,a) \right)^N\!$,
where $\,\vec r_j= \vec r_j(s_j,u_j)\,$ are the Cartesian coordinates
of the $\,N$--th particle.

As we have said our main aim in this paper is to estimate the maximum
number of particles which a curved waveguide with given
$\,\gamma,a\,$ can bind, \ie, to find conditions under which the
discrete spectrum of $\,H_N\,$ is empty. To this end, one has to
determine first the bottom of the essential spectrum. In complete
analogy with the usual HVZ theorem \cite{RS}, we find 
\begin{equation} \label{ess spectrum}
\sigma_{\rm ess}(H_N)\,=\, \left\lbrack\, \mu_{N-1}+ \left(\pi \over
2a \right)^2\!,\, \infty \right)\,,
\end{equation}
where $\,\mu_{N-1}:= \inf\, \sigma(H_{N-1})\,$. Obviously,
$$
\inf\,\sigma_{\rm ess}(H_N)\,\le\, \mu_{N-k} +\, k\left(\pi \over
2a \right)^2
$$
holds for $\,k=1,\dots,N\!-\!1\,$, so
\begin{equation} \label{inf ess spectrum}
\inf\,\sigma_{\rm ess}(H_N)\,\le\, N \left(\pi \over 2a \right)^2.
\end{equation}
In a straight tube the two expressions equal each other, while for
$\,\gamma\ne 0\,$ we have a sharp inequality because $\,\mu_1<
\left(\pi \over 2a \right)^2\,$ holds in this case.


\setcounter{equation}{0}
\section{Neutral fermions}

If the particles in question are neutral fermions, one can get a
simple upper bound on the number of bound states using the
one--particle Hamiltonian (\ref{Hamiltonian 1}); it is sufficient to
estimate the dimension of $\,\sigma_{\rm disc}(H_1)\,$ and to employ
the Pauli principle. To this aim, one has to estimate $\,H_1\,$ from
above by an operator with the transverse and longitudinal variables
decoupled; its projections to transverse modes are then
one--dimensional Schr\"odinger operators to which the modified
Birman--Schwinger bound may be applied 
[17-19]. 
In
Ref.~\cite{DE} we used this argument in the situation where $\,a\,$
is small so that only the lowest transverse mode and the leading term
in (\ref{effective potential}) may be taken into account.

A modification to the more general case is straightforward. We
introduce the function 
\begin{equation} \label{potential bound}
\tilde W(s)\,:=\, {\gamma(s)^2\over 4\delta_-^2}\,+\, {a|\gamma''(s)|
\over 2\delta_-^3} \,+\, {5a^2\gamma'(s)^2\over 4\delta_-^4}\,,
\end{equation}
where
\begin{equation} \label{delta}
\delta_{\pm}\,:=\, 1\pm a\|\gamma\|_{\infty}\,,
\end{equation}
which majorizes the effective potential, $\,V(s,u)\le \tilde W(s)\,$.
Furthermore, we set
\begin{equation} \label{potential bound j}
\tilde W_j(s)\,:=\, \max \left\lbrace\, 0,\, \left(\pi\over 2a
\right)^2 (1\!-\!j^2)\, \right\rbrace
\end{equation}
for $\,j=2,3,\dots\,$; in view of the assumptions {\em (ii), (iii)}
only finite number of them is different from zero.

Replacing $\,V\,$ by $\,\tilde W\,$, and $\,(1\!+\!u\gamma)^{-2}$ by
$\,\delta_+^{-2}$, we get an estimating operator with separating
variables, or in other words, a family of shifted one--dimensional
Schr\"odinger operators; we are looking for the number of their
eigenvalues below $\,\inf \sigma_{\rm ess}(H_1)= \left(\pi\over
2a\right)^2$. The mentioned modification of the Birman--Schwinger
bound is based on splitting the rank--one operator corresponding to
the singularity of the resolvent kernel $\,{1\over 2\kappa}\,
e^{-\kappa|s-s'|}\,$ at $\,\kappa=0\,$ and applying a Hilbert-Schmidt
estimate to the rest.  In analogy with Refs.~\cite{Se,Kl,Ne} 
we employ this trick for the lowest--mode component of the estimating
operator, while for the higher modes we use the full resolvent at the
values $\,\kappa_j:= \left(\pi\over 2a\right) \sqrt{j^2\!-1}\,$.  In
this way we arrive at the following conclusion: 
\begin{proposition}
{\rm The number $\,N\,$ of neutral particles of half--integer spin
$\,S\,$ which a curved quantum waveguide can bind satisfies the
inequality}
\begin{eqnarray} \label{neutral Pauli}
N\,\le\, (2S\!+\!1) \bigg\lbrace\,
1 &\!+\!& \delta_+^2\, {\int_{\R^2} \tilde W(s) |s-t|\, \tilde W(t)
\,ds\,dt \over \int_{\R} \tilde W(s)\,ds} \nonumber \\ \nonumber \\  
&\!+\!& \sum_{j=2}^{\infty}\, {a\delta_+^2\over \pi\sqrt{j^2\!-1}}\:
\int_{\R} \tilde W_j(s)\,ds \,\bigg\rbrace\,.
\end{eqnarray}
\end{proposition}
\begin{remarks}
{\rm (a) As we have said, the number of nonzero term in the last
sum is finite. More exactly, the index $\,j\,$ runs up to the
entire part of $\,\sqrt{1+\,\left(2a\over\pi \right)^2 \|\tilde
W\|_{\infty}}\,$; hence if $\,a\,$ is small enough this term is 
missing at all. \\
(b) The assumption {\em (iv)} is not needed here. It is sufficient,
\eg, that the functions $\,\gamma, \gamma'$, and $\,|\gamma''|^{1/2}$
decay as $\,|s|^{-1-\eps}$ as $\,|s|\to\infty\,$.}
\end{remarks}

\setcounter{equation}{0}
\section{Main result: $N$ charged particles}

We have said in the introduction that the present study is motivated
mainly by the need to describe electrons in curved quantum wires.
Unfortunately, the above simple estimate have no straightforward
consequences for the situation when the particles are charged. While
the electrostatic repulsion adds a positive term to the Hamiltonian
(\ref{Hamiltonian N}), it may move at the same time the bottom
of the essential spectrum since the energies of the bound ``clusters"
are, of course, sensitive to the interaction change.

We need therefore another approach which would allow to take the
repulsion term in (\ref{Hamiltonian N}) into account. An inspiration
can be found in analysis of atomic $\,N$--body Hamiltonians. To
formulate the result we need some notation. Given a positive
$\,\beta\,$ we denote by $\,\{\lambda_m\}_{m=1}^{\infty}$ the ordered
sequence of eigenvalues of Dirichlet Laplacian at the rectangle
\begin{equation} \label{rectangle}
R_{\beta}:=\, \left\lbrack -\,{3\over 2} \beta\delta_+,\, {3\over 2}
\beta\delta_+ \,\right\rbrack\, \times\, [-a,a]\,,
\end{equation}
and set
\begin{equation} \label{T(N)}
T_{\beta}(N)\,:=\, \left\lbrace\, \begin{array}{lll} 2 \sum_{m=1}^n
\lambda_m \quad & \dots & \quad N=2n \\ \\ 2 \sum_{m=1}^n\lambda_m+
\lambda_{n+1} \quad & \dots & \quad N=2n+1 
\end{array} \right.
\end{equation}
We have in mind here electrons and assume that the spin is $\,{1\over
2}\,$, otherwise $\,T_{\beta}(N)\,$ has to be replaced by the sum of
the first $\,N\,$ eigenvalues of $\,2S\!+\!1\,$ identical copies of
the Laplacian. Now we are able state our main result:
\begin{theorem}
{\rm Assume {\em (i)--(iv).} $\,\sigma_{\rm disc}\left(
H_N(\gamma,a,e) \right)= \emptyset\,$ for $\,N\ge2\,$ if the
condition   
\begin{equation} \label{absence}
T_{\beta}(N)\,+\, {e^2\over 2\beta\sqrt{7}}\, N(N\!-\!1) \,\ge\,
\|\tilde W\|_{\infty}N \,+\, \left(\pi\over 2a\right)^2 N \,+\,
{e^2\over 18\beta\sqrt{2}}
\end{equation}
is valid for some $\,\beta\ge \max\{\,2b, 596\,e^{-2}\}\,$.}
\end{theorem}
{\em Proof:} We use a variational argument which relies on a suitable
decomposition of the configuration space. Consider a pair of smooth
functions $\,v,\,g\,$ from $\,\R_+$ to $\,[0,1]\,$ such that
\begin{equation} \label{v}
v(t)\,=\, \left\lbrace\, \begin{array}{lll} 0 & \quad \dots \quad &
t\le 1 \\ \\ 1 & \quad \dots \quad & t\ge {3\over 2} \end{array}
\right. 
\end{equation}
and
\begin{equation} \label{vg decomposition}
v(t)^2\!+g(t)^2=\,1\,.
\end{equation}
Elements of the configuration space are $\,(s,u)\,$ with $\,s=\{s_1,
\dots,s_N\}\,$ and $\,u=\{u_1,\dots,u_N\}\,$. We denote
$\,\|s\|_{\infty}:= \max\{s_1,\dots,s_N\}\,$ and employ the
functions 
$$
s\,\mapsto \,v(\|s\|_{\infty}\beta^{-1}),\;
g(\|s\|_{\infty}\beta^{-1})\,,
$$
where $\,\beta>2b>a\,$ is a parameter to be specified later. By abuse
of notation, we use the symbols $\,v,\,g\,$ again both for these
functions and the corresponding operators of multiplication. It is
straightforward to evaluate  $\,\left([H_N,v]\psi,v\psi\right)\,$ and
the analogous expression with $\,v\,$ replaced by $\,g\,$ for a
vector $\,\psi\in D(H_N)\,$; in both cases it is only the
longitudinal kinetic part in (\ref{Hamiltonian N}) which contributes.
This yields the identity
\begin{eqnarray*}
(H_N\psi,\psi) &\!=\!& (H_Nv\psi,v\psi)\,+\, (H_Ng\psi,g\psi) \\ 
&\!+\!& \sum_{j=1}^N \left\lbrace\, \left\|(1\!+\!u_j\gamma_j)^{-1}
v_j\psi \right\|^2+\, \left\|(1\!+\!u_j\gamma_j)^{-1} g_j\psi
\right\|^2 \,\right\rbrace\,,
\end{eqnarray*}
where we have used the shorthands $\,v_j:= {\partial v\over\partial
s_j}\,$, $\,g_j:= {\partial g\over\partial s_j}\,$, and $\,\gamma_j:=
\gamma(s_j)\,$. Notice further that the factors
$\,(1\!+\!u_j\gamma_j)^{-1}\,$ may be neglected, because $\,v_j\,
g_j\,$ are nonzero only if $\,s_j\ge \beta>2b\,$ in which case
$\,\gamma_j=0\,$. Furthermore, with the exception of the hyperplanes
where two or more coordinates coincide (which is a zero measure set)
the norm $\,\|s\|_{\infty}\,$ coincides with just one of the
coordinates $\,s_1,\dots,s_n\,$, and therefore
\begin{equation} \label{vg error}
\sum_{j=1}^N \left\lbrace\, \left\| v_j\psi \right\|^2+\, \left\|
g_j\psi \right\|^2 \,\right\rbrace\,\le\, \|\psi\|^2 \max_{1\le j\le
N} \left\lbrace\, \left\| v_j \right\|^2_{\infty}+\, \left\| g_j
\right\|^2_{\infty} \,\right\rbrace \,\le\, \beta^{-2}C_0
\|\psi\|^2\,, 
\end{equation}
where $\,C_0:=\|v'\|_{\infty}^2\!+ \|g'\|_{\infty}^2\,$. We arrive at
the estimate
\begin{equation} \label{decomposition 1}
(H_N\psi,\psi)\,\ge\, L_1[v\psi]+ L_1[g\psi]
\end{equation}
with
\begin{equation} \label{L_1}
L_1[\phi]\,:=\, (H_N\phi,\phi)\,-\, {C_0\over\beta^2}\,
\|\phi\|^2_{\NN_\beta}\,,
\end{equation}
where the last index symbolizes the norm of the vector $\,\phi\,$
restricted to the subset $\,\NN_\beta:= \left\lbrace\,s:\;
\beta\le\|s\|_{\infty}\le\, {3\beta \over 2}\,\right\rbrace\,$ of the
configuration space.

Next one has to estimate separately the contributions from the inner
and outer parts. Let us begin with the exterior. We introduce the
following functions:
\begin{eqnarray*}
f_1(s) &\!=\!& v\left(2s_1\|s\|^{-1}_{\infty}\right)\,, \\
f_j(s) &\!=\!& v\left(2s_j\|s\|^{-1}_{\infty}\right)
\prod_{n=1}^{j-1} g\left(2s_n\|s\|^{-1}_{\infty}\right)\,, \qquad
j=2,\dots,N\!-\!1\, \\
f_N(s) &\!=\!& \prod_{n=1}^{N-1}
g\left(2s_n\|s\|^{-1}_{\infty}\right)\,.
\end{eqnarray*}
It is clear from the construction that
\begin{equation} \label{decomposition 2}
\sum_{j=1}^N f_j(s)^2= \,1\,.
\end{equation}
Moreover, the functions
$$
s_j\,\mapsto\, v(2s_j\|s\|_{\infty}),\, g(2s_j\|s\|_{\infty})
$$
have a non--zero derivative only if $\,|s_j|\ge\, {1\over
2}\|s\|_{\infty}^{-1}$. Hence on the support of $\,s\mapsto
v(\|s\|_{\infty}\beta^{-1})$ the derivative is non--zero if
$\,|s_j|\ge \,{1\over 2}\beta>b\,$. In other words, the function
$\,s\mapsto f_j(s)^2 v(\|s\|_{\infty}\beta^{-1})$ has zero derivative
in all the parts of the configuration space where at least one of the
electrons dwells in the curved part of the waveguide. Commuting the
(longitudinal kinetic part of) $\,H_N\,$ with $\,f_j\,$, we get in
the same way as above the identity
\begin{equation} \label{decomposition 3}
L_1[v\psi]\,=\, \sum_{j=1}^N \left\lbrace\, L_1[f_jv\psi] -
\|(\nabla_sf_j)v\psi\|^2 \,\right\rbrace\,,
\end{equation}
where $\,\nabla_s:= \left(\partial_{s_1},\dots, \partial_{s_1}
\right)\,$. Next we need a pointwise upper bound on $\,\sum_{j=1}^N
(\nabla_sf_j)^2$: denoting $\,\sigma_j:= 2s_j\|s\|_{\infty}$, we can
write 
\begin{eqnarray*}
\sum_{j=1}^N |(\nabla_sf_j)(s)|^2 &\!=\!& {4\over\|s\|_{\infty}^2}\,
\bigg\lbrace\, v'(\sigma_1)^2 \\ \\
&& +g'(\sigma_1)^2 v(\sigma_2)^2 +g(\sigma_1)^2 v'(\sigma_2)^2 
+ \cdots \\
&& +g'(\sigma_1)^2 g(\sigma_2)^2\!\dots g(\sigma_N)^2 +\cdots+
g(\sigma_1)^2 \!\dots g(\sigma_{N-1})^2 g'(\sigma_N)^2\,
\bigg\rbrace\,, 
\end{eqnarray*}
which gives after a partial resummation
\begin{eqnarray*}
&\!=\!& {4\over\|s\|_{\infty}^2}\, \bigg\lbrace\, v'(\sigma_1)^2 
+ g'(\sigma_1)^2 +g(\sigma_1)^2 g'(\sigma_2)^2 + \cdots \\
&& +g(\sigma_1)^2 \!\dots g(\sigma_{N-1})^2 g'(\sigma_N)^2\,
\bigg\rbrace \\ \\
&\!\le\!& {4\over\|s\|_{\infty}^2}\, \left\lbrace\, v'(\sigma_1)^2
+\sum_{j=1}^N g'(\sigma_j)^2  \,\right\rbrace\,\le\,
{4NC_0\over\|s\|_{\infty}^2} \;;
\end{eqnarray*}
recall that $\,C_0:=\|v'\|_{\infty}^2\!+ \|g'\|_{\infty}^2\,$.
Consequently, 
\begin{eqnarray} \label{L_1 decomposition}
L_1[v\psi] &\!\ge \!& \sum_{j=1}^N L_1[f_j v\psi] \,-\, 4NC_0 \left\|
v\psi \|s\|_{\infty}^{-1} \right\|^2 \nonumber \\ \nonumber \\
&\! =\!& \sum_{j=1}^N\, \left\lbrace\, L_1[f_j v\psi] -4NC_0 \left\|
f_j v\psi \|s\|_{\infty}^{-1} \right\|^2\, \right\rbrace  \nonumber \\
\nonumber \\ 
&\! =\!& \sum_{j=1}^N\, L_2[f_j v\psi]\,,
\end{eqnarray}
where
\begin{equation} \label{L_2}
L_2[\phi]\,:=\, L_1[\phi] \,-\,4NC_0 \left\| \phi \|s\|_{\infty}^{-1}
\right\|^2 \,.
\end{equation}
Hence we have to find a lower bound to $\,L_2(\psi_j)\,$ with
$\,\psi_j:= f_jv\psi)\,$. Since $\,s_j\ge\, {1\over 2}\|s\|_{\infty}
\ge\, {1\over 2}\beta >b\,$ holds on the support of $\,\psi_j\,$, we have
$\,V(s_j,u_j)=0\,$ there. This allows us to write
$$
(H_N\psi_j,\psi_j)\,=\, (H_{N-1}\psi_j,\psi_j)+\, \left\|
\partial_{s_j} \psi_j \right\|^2+\, \left\| \partial_{u_j} \psi_j
\right\|^2+ e^2\, \sum_{j\ne l=1}^N \left( |\vec r_j\!- \vec
r_l|^{-1} \psi_j, \psi_j \right)\,, 
$$
where $\,H_{N-1}$ refers to the system with the $\,j$--th electron
excluded, and therefore
$$
(H_N\psi_j,\psi_j)\,\ge\, \left(\, \mu_{N-1}+\, \left( \pi\over 2a
\right)^2 \right) \|\psi_j\|^2 +\, e^2\, \sum_{j\ne l=1}^N \left(
|\vec r_j\!- \vec  r_l|^{-1} \psi_j, \psi_j \right)\,.
$$
Since $\,|\vec r_j\!- \vec r_l|\le \sqrt{ (s_j\!-\!s_l)^2\! +4a^2}\,
\le\, 2\, \sqrt{\|s\|_{\infty}^2\!+a^2}\,$, we have 
$$
(H_N\psi_j,\psi_j)\,\ge\, \left(\, \mu_{N-1}+\, \left( \pi\over 2a
\right)^2 \right) \|\psi_j\|^2 +\, {e^2(N\!-\!1)\over 2}\, \left(
(\|s\|^2\!+a^2)^{-1/2} \psi_j,\psi_j \right)\,.
$$
The sought lower bound then follows from (\ref{L_2}) and (\ref{L_1}):
\begin{eqnarray*}
L_2[\psi_j] &\! \ge \!& \left(\, \mu_{N-1}+\, \left( \pi\over 2a
\right)^2 \right) \|\psi_j\|^2 -\,4NC_0 \left\| \psi_j
\|s\|_{\infty}^{-1} \right\|^2\\ \\
&\! -\!& C_0\beta^{-2} \|\psi_j\|^2_{\NN_\beta} +\,
{e^2(N\!-\!1)\over 2}\, \left( (\|s\|^2\!+a^2)^{-1/2} \psi_j,\psi_j
\right)\;; 
\end{eqnarray*}
recall that $\,\NN_\beta:= \left\lbrace\,s:\;
\beta\le\|s\|_{\infty}\le\, {3\beta \over 2}\,\right\rbrace\,$. The
second and the third term at the \rhs can be combined using
$$
4NC_0 \left\| \psi_j \|s\|_{\infty}^{-1}\right\|^2 +\,C_0 \beta^{-2}
\|\psi_j\|^2_{\NN_\beta} \,\le\, (4N\!+\!1)C_0 \left\| \psi_j
\|s\|_{\infty}^{-1}\right\|^2\,. 
$$
Furthermore, $\,\|s\|_{\infty}\ge \beta >2b>a\,$ yields
$\,(\|s\|^2\!+a^2)^{1/2} \le \sqrt{2}\, \|s\|_{\infty}$ and
\begin{eqnarray} \label{L_2 bound}
L_2[\psi_j] &\!\ge\!& \left(\, \mu_{N-1}+\, \left( \pi\over 2a
\right)^2 \right) \|\psi_j\|^2 \nonumber \\ \nonumber \\ 
&\!+\!&  \left(\, {e^2(N\!-\!1)\over 2\sqrt{2}} \,-\,
{C_0(4N\!+\!1)\over \beta}\,\right)\, \left\| \psi_j
\|s\|_{\infty}^{-1}\right\|^2 \,. 
\end{eqnarray}
We are interested in the situation when the second term at the \rhs
is positive. This is achieved if
$$
{e^2(N\!-\!1)\over 2\sqrt{2}} \,>\, {C_0(4N\!+\!1)\over \beta}
$$
which is ensured if we choose $\,\beta\,$ in such a way that
\begin{equation} \label{beta bound}
\beta\,>\, {18\sqrt{2} C_0\over e^2}\;;
\end{equation}
recall that $\,N\ge 2\,$. Owing to the identity (\ref{L_1
decomposition}) we then have
\begin{equation} \label{outer bound}
L_1[v\psi]\,\ge\, \left(\, \mu_{N-1}+\, \left( \pi\over 2a
\right)^2 \right) \|v\psi\|^2\,,
\end{equation}
which means in view of (\ref{ess spectrum}) that the external part of
$\,\psi\,$ does not contribute to the discrete spectrum.

Let us turn now to the inner part. The corresponding quadratic form
in the decomposition (\ref{decomposition 1}) can be estimated with
the help of (\ref{Hamiltonian N}) and (\ref{L_1}) by
\begin{eqnarray} \label{inner bound}
L_1[g\psi] &\! \ge \!& \delta_+^{-2} \|\nabla_s g\psi\|^2+\,
\|\nabla_u g\psi\|^2 +\, \sum_{j=1}^N \left( V(s_j,u_j)g\psi, g\psi
\right) \nonumber \\ \nonumber \\
&\! +\!& e^2 \sum_{1\le<k\le N} \left( \|\vec r_j\!-\vec r_k\|^{-1}
g\psi, g\psi \right) \,-\, {C_0\over \beta^2}\, \|g\psi\|^2\;;
\end{eqnarray}
recall that $\, \delta_+:= 1\!+\!a\|\gamma\|_{\infty}$. Using the
function $\,\tilde W\,$ defined by (\ref{potential bound}) we find
$\,|V(s_j,u_j)| \le \tilde W(s_j)\,$, so
$$
\max \left\lbrace\, V(s,u)\,:\; (s,u)\in \R\times [-a,a]\,
\right\rbrace\, \le\, \|\tilde W\|_{\infty}\,.
$$
Consequently, the curvature--induced potential term can be estimated
by 
$$
\sum_{j=1}^N\, \left( V(s_j,u_j) g\psi, g\psi \right)\,\le\, \|\tilde
W\|_{\infty}\, N \|g\psi\|^2\,.
$$
Furthermore, on the support of $\,g\,$ we have
$$
|\vec r_j\!- \vec r_k| \le\, 2\, \sqrt{\|s\|_{\infty}^2\!+a^2}\,
\le\, \sqrt{3\beta^2\!+4a^2}\,, 
$$
because $\,\|s\|_{\infty}\le\, {3\over 2}\beta\,$ holds there. At the
same time, $\,\beta>2b>a\,$, so we arrive at the estimate
$$
|\vec r_j\!- \vec r_k| \le\, \sqrt{7} \beta\,,
$$
which yields
$$
\sum_{1\le<k\le N} \left( \|\vec r_j\!-\vec r_k\|^{-1}
g\psi, g\psi \right)\,\ge\, {N(N\!-\!1)\over 2\beta\sqrt{7}}\,
\|g\psi\|^2\,.
$$
Now we can combine the above estimates with the inequality
$\,{C_0\over \beta}\,<\, {e^2\over 18\sqrt{2}}\,$ which follows from
(\ref{beta bound}) to get the bound
\begin{eqnarray} \label{inner bound 2}
L_1[g\psi] &\! \ge \!& \delta_+^{-2} \|\nabla_s g\psi\|^2+\,
\|\nabla_u g\psi\|^2  \nonumber \\ \nonumber \\
&\! +\!& \left\lbrack\, -N\|\tilde W\|_{\infty} +\, {e^2
N(N\!-\!1)\over 2\beta \sqrt{7}}\,-\, {e^2 \over 18\beta \sqrt{2}}\,
\right\rbrack \, \|g\psi\|^2\,.
\end{eqnarray}
Now we can put the above results together. In view of the inequality
(\ref{outer bound}) and of (\ref{inf ess spectrum}), the last bound
tells us that $\,H_N\,$ has no discrete spectrum for $\,N\ge 2\,$
provided 
\begin{eqnarray} \label{discrete absence}
\delta_+^{-2} \|\nabla_s g\psi\|^2 &\!+\!& \|\nabla_u g\psi\|^2+
\,\bigg\lbrack\, {e^2
N(N\!-\!1)\over 2\beta \sqrt{7}}\,-\, {e^2 \over 18\beta \sqrt{2}}
\nonumber \\  \nonumber \\ &\!-\!&
N\|\tilde W\|_{\infty} -\, N \left(\pi \over 2a\right)^2
\bigg\rbrack \, \|g\psi\|^2 \,\ge\, 0
\end{eqnarray}
for some $\,\beta\,$ which satisfies the condition
\begin{equation} \label{beta bound 2}
\beta\,\ge\, \max\, \left\lbrace\, 2b,\, {18\sqrt{2} C_0\over e^2}\,
\right\rbrace\,. 
\end{equation}
The first two terms in (\ref{discrete absence}) are nothing else than
the quadratic form of the $\,2N$--dimensional Laplacian on
$\,R_{\beta}^N$ --- \cf (\ref{rectangle}). By Pauli principle each
eigenvalue may appear only twice, thus one has to take the orthogonal
sum of two copies of the Laplacian on $\,R_{\beta}$ and to summ the
first $\,N\,$ eigenvalues of such an operator. This is exactly the
quantity which we have called $\,T_{\beta}(N)\,$. 

To finish the proof, it remains to estimate $\,C_0\,$ which appears
in the conditions (\ref{beta bound}) and (\ref{beta bound 2}). We
will not attempt an optimal bound and put simply  
$$
v(\xi)\,:=\, \sin \left(4\pi\xi^2(1\!-\!2\xi^2)\right)
$$
for $\,t\!-\!1=:\xi\in \left(0,{1\over 2}\right)\,$, then 
$$
v'(\xi)^2+ g'(\xi)^2=\, (8\pi)^2\xi^2(1-4\xi^2)^2
$$
has the maximum value $\,2\sqrt{2}(8\pi)^2/3\approx 595.5\,$. \quad
\QED 


\setcounter{equation}{0}
\section{Conclusions}

Since the present study is rather a foray into an unchartered
territory, the result is naturally far from optimal. Let us add a few
remarks. First of all, it is clear that the overall size of the
curved region affects substantially the number of particles which the
waveguide can bind. We know that {\em any} curved tube has a
one--particle bound state \cite{DE,GJ}, hence a tube with $\,N\,$
slight bends which very far from each other (so far that the
repulsion is much smaller that the gap between the bound state energy
and the continuum) can certainly bind $\,N\,$ particles for $\,N\,$
arbitrarily large.

The method we use is borrowed from atomic physics where it yields
bounds on ionization of an atom. Of course, there are differences.
The binding is due to the curved hard wall of the waveguide rather
than by the electrostatic attraction to the nucleus, and the spectrum
of our one--particle operator (\ref{Hamiltonian 1}) is finite.
Consequently, there is a maximum number of particles which a given
curved tube can bind as long as the particles are fermions. Bosons
can occupy naturally a single state, and the idea of a {\em Bose
condensate} of neutral spin--zero atoms in a curved hollow optical
fiber is rather appealing.

On the other hand, a non--zero particle charge changes the picture,
and even the number of bosons bind by a curved tube is limited:
notice that the condition (\ref{absence}) is satisfied for large
enough $\,N\,$ without respect to the Pauli--principle term
$\,T_{\beta}(N)\,$. Of course, the fermionic nature reduces the
maximum number $\,N\,$ further, since $\,T_{\beta}(N)\,$ growth for
large $\,N\,$ is between $\,o(N^3)\,$ in the limit $\,a\to 0\,$ and
$\,o(N^2)\,$ for $\,2b\sim a\,$. At the same time, the maximum number
also depends on the value of the charge. Since $\,{1\over\sqrt{7}}
\,-\, {1\over 18\sqrt{2}} \,>0\,$ and the remaining terms in
(\ref{absence}) are independent od $\,e\,$ we see that $\,\sigma_{\rm
disc}(H_N)= \emptyset\,$ for any $\,N\ge 2\,$ provided $\,e\,$ is
large enough. Thus our result confirms the natural expectation that
for a given curved tube and sufficiently charged particles just
one--particle bound states can survive.

We have not addressed in this paper the question about the minimum
number of particles which a curved quantum waveguide can bind. The
gap between the trivial result which follows from the one--particle
theory \cite{DE,ES,GJ} and the condition (\ref{absence}) leaves a lot
of space for improvements. Moreover, it is a natural question whether
strongly curved tubes which can bind many particles allow for some
semiclassical description analogous to the case of the quantum dots
\cite{Y}. This is a task for a future work.


\section*{Acknowledgments}

The research has been partially supported by the Grants
No.~202--0218, GACR, and ME099, Ministry of Education of the Czech
Republic. 

\end{document}